\newcommand{\calC}{\mathcal{C}}

\newcommand{\calE}{\mathcal{E}}
\newcommand{\calF}{\mathcal{F}}
\newcommand{\calG}{\mathcal{G}}

\newcommand{\calL}{\mathcal{L}}

\newcommand{\calT}{\mathcal{T}}
\newcommand{\calU}{\mathcal{U}}
\newcommand{\calO}{\mathcal{O}}

\newcommand{\calZ}{\mathcal{Z}}

\newcommand{\bbD}{\mathbb{D}}

\newcommand{\bbX}{\mathbb{X}}

\newcommand{\bbZ}{\mathbb{Z}}

\newcommand{\bbF}{\mathbb{F}}
\newcommand{\bbC}{\mathbb{C}}
\newcommand{\bbH}{\mathbb{H}}

\newcommand{\bbP}{\mathbb{P}}
\newcommand{\bbQ}{\mathbb{Q}}
\newcommand{\bbR}{\mathbb{R}}
\newcommand{\bbV}{\mathbb{V}}
\newcommand{\Z}{\mathbb{Z}}

\newcommand{\bfs}{\mathbf{s}}

\newcommand{\bfF}{\mathbf{F}}
\newcommand{\bfE}{\mathbf{E}}
\newcommand{\la}{\langle}
\newcommand{\ra}{\rangle}

\newcommand{\SL}{\textrm{SL}}
\newcommand{\Aut}{\text{Aut}}
\newcommand{\Hom}{\text{Hom}}
\newcommand{\Char}{\text{Char}}

\newcommand{\GL}{\text{GL}}

\newcommand{\Pic}{\text{Pic}}
\newcommand{\Ker}{\text{Ker}}

\newcommand{\Spec}{\text{Spec}}

\newcommand{\PSL}{\text{PSL}}

\newcommand{\rank}{\text{rank}}

\newcommand{\Coker}{\textup{Coker}}

\newcommand{\Tors}{\textup{Tors}}
\newcommand{\SU}{\textup{SU}}
\newcommand{\U}{\textup{U}}
\newcommand{\PSU}{\textup{PSU}}

\newcommand{\tga}{\Pi}
\newcommand{\ga}{\Gamma}

\newcommand{\beq}{\begin{equation}}
\newcommand{\eeq}{\end{equation}}
\newcommand{\ab}{\textup{ab}}

\newcommand{\vdeg}{\textup{vdeg}}
\newcommand{\Hilb}{\textup{Hilb}}
\newcommand{\lcm}{\textup{lcm}}
\newcommand{\coh}{\textup{coh}}

\documentclass[11pt]{amsart}
\usepackage{graphicx}
\usepackage{amscd, amsmath, amsthm, amssymb,xy, xypic}

\xyoption{arrow}

\xyoption{matrix}
\input xy
\xyoption{all}

\newtheorem{theorem}{Theorem}[section]
\newtheorem{lemma}[theorem]{Lemma}
\newtheorem{proposition}[theorem]{Proposition}
\newtheorem{conjecture}[theorem]{Conjecture}

\theoremstyle{definition}     
\newtheorem{definition}[theorem]{Definition}

\newtheorem{example}[theorem]{Example}

\theoremstyle{remark}
\newtheorem{remark}[theorem]{Remark}

\numberwithin{equation}{section}
\title[McKay correspondence]{McKay's correspondence for cocompact discrete subgroups of $\SU(1,1)$}
\author{Igor V. Dolgachev}
\address{Department of Mathematics, University of Michigan, 525 E. University Av., Ann Arbor, Mi, 49109}
\email{idolga@umich.edu}
\thanks{The  author was supported in part by NSF grant 0245203.}

\dedicatory{To John McKay}

\begin{document}
\maketitle

\begin{abstract}
The classical McKay correspondence establishes an explicit link from the representation theory  of a finite subgroup $\Pi$ of $\SU(2)$ and the geometry of the minimal resolution of the affine surface $V = \bbC^2/\Pi$. In this paper we discuss a possible generalization of the McKay correspondence to the case when $\Pi$  is replaced with a discrete cocompact subgroup of the universal cover of $\SU(1,1)$ such that its image $\Gamma$ in $\PSU(1,1)$ is a fuchsian group of signature $(0,e_1,\ldots,e_n)$. We establish a correspondence between  a certain class of finite-dimensional unitary representations of $\Pi$ and vector bundles on an open algebraic surface with  trivial canonical class canonically associated to $\Gamma$. 
\end{abstract}

\section{Introduction} It has been known since the work of P. Du Val in the thirties that  Coxeter-Dynkin diagrams  of type $ADE$ are in  bijective correspondence with the conjugacy classes of finite subgroups $\Pi$ of $\SU(2)$ in such a way that the intersection graph of a minimal resolution of $\bbC^2/\Pi$ is the diagram corresponding to the group $\Pi$. In the early eighties John McKay added more to this mysterious connection by introducing a certain graph attached to any finite group. When the group is equal to a binary polyhedral group $\Pi$ the graph coincides with  the affine extension of the  Dynkin diagram attached  to $\Pi$ \cite{McK1},\cite{McK2}.  The vertices of the graph correspond to irreducible representations of $\Pi$ where the extended vertex corresponds to the trivial representation. The first geometric explanations of the McKay correspondence were  given independently by G. Gonzalez-Sprinberg and J.-L. Verdier \cite{GV} and H. Kn\"orrer \cite{Kn}. Other, more algebraic, interpretations were given later  by B. Kostant \cite{Ko}, T. Springer \cite{Sp}, R. Steinberg \cite{St},\cite{St2},  M. Artin and J.-L. Verdier \cite{AV}, H. Esnault and H. Kn\"orrer \cite{EK}.

Modern development reveals a  more general context of the correspondence. The current slogan is that the McKay correspondence establishes an isomorphism of the Grothendieck group $K_G(X)$ of $G$-equivariant coherent sheaves on an algebraic variety $X$ on which a finite group $G$ acts and the Grothendieck group $K(Y)$ of coherent sheaves on a crepant resolution of the quotient $X/G$ (when it exists), or more generally, an equivalence of the corresponding derived categories. For example, when $X= \bbC^n, n = 2,3$ a crepant resolution exists if and only if $G$ is a subgroup of $\SU(n)$ and an equivalence of the categories was established by M. Kapranov and Vesselot ($n = 2$) and T. Bridgeland, A. King and M. Reid ($n= 3)$. We refer for all this and much more to an excellent survey of M. Reid \cite{Reid}.

In this paper we propose a generalization of the McKay correspondence to a certain class of discrete infinite groups $\Pi$. It is based on the following observations. For any finite subgroup $\Pi$ of $\SU(2)$, the affine surface $\bbC^2/\Pi$ is isomorphic to the affine spectrum of the algebra of $\Gamma$-invariant sections of the tangent bundle on a simply connected Riemann surface $\bbP^1$, where $\Gamma$ is the image of $\Pi$ in $\PSU(2)$. Replacing $\bbP^1$ with the  unit disk $\bbD$ and the tangent bundle with the cotangent bundle, we can consider the algebra of automorphic forms $A(\Gamma)$ with respect to any cocompact discrete subgroup $\Gamma$ of $\PSU(1,1)$ of some signature $(g;e_1,\ldots,e_n)$. It is a finitely generated graded algebra and its affine spectrum  is an affine surface $V_\ga$ with an isolated singular point $0\in V_\Gamma$ corresponding to the unique maximal homogeneous ideal. The fundamental group of $V_\Gamma^* = V_\ga\setminus \{0\}$ is isomorphic to a central extension $\Pi$ of $\Gamma$ equal to the pre-image of $\Gamma$ in the universal cover of $\PSU(1,1)$. This suggests to consider   linear representations of $\Pi$ as the first side of the McKay correspondence. In fact,  to keep the analogy closer,  we have to restrict ourselves to  finite-dimensional unitary linear representations of $\Pi$.  Each such representation  is isomorphic to the pull-back of a  unitary representation of the fundamental group of a Deligne-Mumford stacky hyperbolic curve (see \cite{CS}).  Note that a linear irreducible representation of a binary polyhedral group is  isomorphic to the pullback of a unitary  representation of a spherical DM-curve. Any unitary representation of $\Pi$ is isomorphic to the direct sum of irreducible ones. We restrict ourselves to admissible representations. i.e.  representations isomorphic to the direct sum of irreducible unitary representations whose image intersects $\SL(n,\bbC)$ in a Zariski dense subgroup. In particular, we want to avoid representations that are factored through a subgroup of finite index of $\Pi$.

To get the second side of the McKay correspondence we replace a minimal resolution of $V_\Gamma$ (never crepant)  with an open  subset $\widetilde{V}_\Gamma$ of  a unique minimal smooth normal crossing  compactification of  $V_\Gamma$ preserving the $\bbC^*$-action defined by the grading of $A(\ga)$. Its complement is the exceptional curve $\bfE_0$ of the minimal resolution  of the singular point $0$ of $V_\Gamma$.  It has been known to the author for 30 years (as well as to some other experts) that,  in the case when $\bbD/\Gamma \cong \bbP^1$, the curve $\bfE_\infty = \widetilde{V}_\Gamma\setminus V_\Gamma^*$ is equal to the union of smooth rational curves with self-intersection $-2$ and its dual graph is a star-shaped tree with $n$ arms of length $e_1,\ldots,e_n$ (counting  the center). In particular, the canonical class of $\widetilde{V}_\Gamma$ is zero. So the surface $ \widetilde{V}_\Gamma$ seems to be a correct substitute for the minimal resolution of singularities of the surface $V_\Gamma$ in the case of finite $\Gamma$.

Let $\rho$ be an admissible irreducible representation of $\Pi$ and $m$ be the order of the image of the center of $\Pi$ (the \emph{level} of $\rho$). We choose a surface subgroup $\Gamma'$ of $\Gamma$ with $C = \bbD/\Gamma'$ of genus $g$ and $m| 2g-2$. The restriction  of $\rho$ to the pre-image $\Pi'$ of $\Gamma'$ in $\Pi$  defines an irreducible  local unitary system on the curve $C$. By a theorem of Narasimhan and Seshadri \cite{NS}, it correspondence to a stable vector bundle $E$ on $C$. We show that $E$ admits a linearization with respect to some central extension $m.G$  of the group $G = \ga/\ga'$. We take the pull-back $\pi^*(E)$ of $E$ to the ruled surface $\pi:X\to C$ defined by a $m$-root of the tangent bundle and then restrict $\pi^*(E)$ to the complement  $X^0$ of the exceptional section $S$. The quotient 
$X^0/m.G$ has only cyclic quotient singularities of types $A_{e_i}$.  The usual McKay correspondence assigns to $\pi^*(E)|X^0$ a unique vector bundle on $\widetilde{V}_\ga$. This is our generalized McKay correspondence.

 In the case when $\Gamma$ is of signature $(0;2,3,7)$ one can choose $C$ to be the Klein quartic curve of genus $3$ and $G = \PSL(2,\bbF_7)$. We showed in \cite{Do1} that there are two stable $\SL(2,\bbF_7)$-linearized bundles of rank 2 and also two of rank 3. Thus there are two admissible representations of $\Pi$ of dimension 2 of level 2 and two of dimension 3 of level 1. It is clear that we can construct admissible representations of any dimensions, for example in the Klein case, by considering symmetric powers of the two-dimensional representations. We conjecture that there exists a certain set of ``fundamental admissible representations'' bijectively corresponding to the irreducible components of the curve $\bfE_\infty$ making an analogy with fundamental weights of a simple infinite-dimensional hyperbolic Lie algebra with the Dynkin diagram equal to the dual graph of the curve $\bfE_\infty$. For example, in the Klein case the diagram is $E_{10}$ so we are looking for 10 fundamental admissible irreducible representations of $\Pi$.

 I am grateful to Daniel Allcock for many useful conversations on the topic of this paper. I thank Brian Jurgelewicz for his help in making the manuscript more reader-friendly.

 \section{Canonical affine surfaces with a good $\bbC^*$-action} \label{sec1}
A familiar construction of ADE-singularities is as follows. Let 
$$h = |z_1|^2+|z_2|^2:\bbC^2\to \bbR$$ be a positive definite Hermitian from, and
$$\calU  = h^{-1}(\bbR_{>0}) = \bbC^2\setminus \{0\}.$$
The projection $\calU\to \bbP^1, (z_1,z_2)\mapsto z_2/z_1$ is surjective. The fibres are isomorphic to $\bbC^*$. Let 
$$\Pi\subset  SU(2)$$
be a finite subgroup acting naturally on $\calU$. Its image $\Gamma$ in $\PSU(2)$ is a polyhedral group, a finite group of automorphisms of $\bbP^1$. The kernel $\Pi\to \Gamma$ is the center $\{\pm 1\}$ of $\SU(2)$ (unless $\Gamma$ is a cyclic group of odd order, in which case the kernel is trivial).  Consider the orbit space $\calU/\Pi.$ Then there exists a normal affine surface $V_\Gamma$ with a singular point $0\in V_\Gamma$ such that 
\begin{equation}\label{quot1}
V_\Gamma^* = V_\Gamma\setminus \{0\} \cong \calU/\Pi\end{equation}
The singular point $(V_\Gamma ,0)$ represents the analytical class of a double rational point, or an ADE-singularity, or a Du Val singularity. 

The surface $V_\Gamma$ is naturally isomorphic to the graded algebra 
\beq\label{ag}
A(\Gamma) = \bigoplus_{k=0}^\infty H^0(\bbP^1,T_{\bbP^1}^{\otimes k})^\Gamma,
\eeq
where $T_{\bbP^1}$ is the  tangent bundle  on $\bbP^1$.
The surface $T_{\bbP^1}^\# = T_{\bbP^1}\setminus \{\text{zero section}\}$ has  universal cover isomorphic to $\calU$ with  Galois group equal to the center of $\Pi$. The group $\Pi$ is isomorphic to the fundamental group of the surface $V_\Gamma^*$.

Now we change the signature of $h$. Let 
$$h =  |z_1|^2-|z_2|^2:\bbC^2\to \bbR$$
$$\calU =  h^{-1}(\bbR_{>0}).$$
The projection $\calU\to \bbP^1, (z_1,z_2)\mapsto z_2/z_1$ has the image equal to the unit disk $\bbD = \{z\in \bbC:|z| < 1\}$. It can be naturally identified with the complement of the zero section of the cotangent bundle $T_{\bbD}^*$. An analog of a  finite subgroup of $\PSU(2)$ in this case is a {\it cocompact discrete subgroup}  
$\ga \subset \PSU(1,1).$ 
Let 
$\Pi$ be its pre-image in the universal cover $\widetilde{\SU}(1,1)$ of $\SU(1,1)$ and $\Pi'$ be the image of $\Pi$ in $\SU(1,1)$. The group $\Pi'$ is a central double extension of $\ga$. 
The group $\Pi'$ acts discretely on $\calU$ and we can take the quotient $\calU/\Pi'$.

\begin{theorem} There exists a unique normal affine surface $V_\Gamma$ with a good $\bbC^*$-action such that 
\begin{equation}\label{quot11}
V_\Gamma^* = V_\Gamma\setminus \{0\} \cong \calU/\Pi'
\end{equation}\eqref{quot1}, where $0$ is the unique closed orbit of the $\bbC^*$ action. 
\end{theorem}

We replace $\bbP^1$ with $\bbD$, another simply-conencted complex manifold of dimension 1.  The algebra \eqref{ag} is now replaced with the graded algebra of automorphic forms with respect to $\Gamma$
\beq\label{ag2}
A(\Gamma) = \bigoplus_{k=0}^\infty H^0(\bbD,T_\bbD^*{}^{\otimes k})^\Gamma.
\eeq
The surface $V_\Gamma$ is defined to be the affine spectrum of $A(\Gamma)$. The group $\Gamma$ acts on $(T_\bbD^*)^\# = T_\bbD^*\setminus \{\text{zero section}\}$ freely and the quotient is isomorphic to $V_\Gamma^*$. The universal cover of 
$V_\Gamma^*$ is equal to the universal cover of $\calU$ and the fundamental group of $V_\Gamma^*$ coincides with the central extension $\Pi$ of $\Gamma$, the pre-image of $\Gamma$ in the universal cover of $\PSU(1,1)$. 

We call the surface $V_\Gamma$ the \emph{canonical affine surface} defined by the group $\Gamma$. The grading of  its ring of regular functions defines a good $\bbC^*$-action, i.e. the action of $\bbC^*$ with a unique closed orbit.

\begin{remark}  It is known that any normal affine surface $V$  with  a good $\bbC^*$-action is obtained in a similar way. One considers a holomorphic line bundle $L$ on a simply connected Riemann surface $P$ admitting an ample  linearization with respect to a cocompact subgroup $\Gamma$ of $\Aut(P)$ acting freely outside the zero section of $L$. Then one proves that the graded algebra 
$$A(\Gamma;L) = \bigoplus_{k=0}^\infty H^0(P,L^{\otimes k})^\Gamma$$
is finitely generated and hence defines a normal affine surface with a good $\bbC^*$-action. Any such surface is obtained from a unique pair $(\Gamma,L)$ \cite{Do1},\cite{Pi}. In the case when $P = \bbP^1$ and $L = T_P$, the algebra $A(P,L)$ is always generated by 3 elements and hence the surface $V_\Gamma$ is isomorphic to a hypersurface in $\bbC^3$. In the case when $P = \bbD$ this could happen only if  some power of $L$ is $\Gamma$-isomorphic to $T_P^*$ \cite{Do2}. There are 31 possible signatures of $\Gamma$ when it really happens (see  \cite{Sch}, \cite{Wa}). Fourteen of them have signature 
$(0;e_1,e_2,e_3)$ and the singular point of $V_\Gamma$ in this case is analytically isomorphic to one of the 14 exceptional unimodal singularities of V. Arnold \cite{Do0} (see  nice accounts of these results in  \cite{Anna}, \cite{Wa2}). 
\end{remark}

Recall that a cocompact discrete subgroup of $\PSU(1,1)$ is described by its signature
$(q;e_1,\ldots,e_n )$, where $q$ is the genus of the Riemann surface $\bbD/\ga$ and $e_i$ are the orders of the non-trivial stabilizer subgroups of $\ga$. The cover $\bbD\to \bbD/\ga$ is a Galois cover ramified over $n$ points 
$p_1,\ldots,p_n$ in $\bbD/\ga$ with ramification indices $e_1,\ldots,e_n $. 
One of the differences between the finite and infinite case is that $\calU$ is simply-connected in the first case and is not simply-connected in the second case. In fact, the group $\SU(2)$ is simply-connected but $\SU(1,1)$ is not. The universal cover $\widetilde{\SU}(1,1)$ can be described  as follows (where I am following a nice exposition from \cite{Anna}). First identify a matrix
$$\begin{pmatrix}z&w\\
\bar{w}&\bar{z}\end{pmatrix}\in \SU(1,1)$$
with an element of the  set 
$X= \{(z,w)\in \bbC^2:|z|-|w|^2 = 1\}$. Then 
$$\widetilde{X} = \{(w,\phi,r)\in \bbC\times \bbR\times \bbR_+:|w|^2 = r^2-1\}$$
with the map $(w,\phi,r) \mapsto (re^{i\phi},w)$ is the universal covering of $X$. Equivalently, we may replace the unit disk $\bbD$ with the upper half-plane $\bbH$, the group $\SU(1,1)$ with $\SL(2,\bbR)$ and realize the universal cover $\widetilde{\SL}(2,\bbR)$ as the subgroup
$$\{(\alpha,g)\in  \calO(\bbH)^*\times \SL(2,\bbR):e^{2\pi i\alpha} = \tfrac{dg}{dz}\}$$
of the semi-direct product $\calO(\bbH)^*\rtimes \SL(2,\bbR)$, where $\SL(2,\bbR)$ acts naturally on $\bbH$ and hence on the group of invertible holomorphic functions $\calO(\bbH)^*$.

Let $\widetilde{\SU}(1,1)$ be the universal cover. It is a 3-dimensional simply-connected Lie group given by the universal central extension
$$1\to \Z \to \widetilde{\SU}(1,1) \to \PSU(1,1) \to 1.$$
The pre-image of $\ga$ in $\widetilde{\SU}(1,1)$ is the group $\tga$ which fits in the central extension
\beq\label{univext}
1\to \Z \to \tga \to \ga\to 1.
\eeq
Recall that 
$$\pi_1(V_\Gamma^*) \cong \tga.$$
It is known that the link space of the singularity $0\in V_\ga$ is a Seifert 3-manifold with invariants 
$(q;(e_1,1),\ldots,(e_n,1))$. The group $\ga$ has a standard set of generators $g_1,\ldots,g_n, a_1,\ldots,a_q,b_1,\ldots,b_q$ with defining relations
$$g_1^{e_1}= \ldots = g_r^{e_n } = g_1\cdots g_n [a_1,b_1]\cdots [a_q,b_q] = 1.$$
The group $\tga$ has generators $c, \tilde{g_1},\ldots,\tilde{g_n}, \tilde{a}_i, \tilde{b}_i, i = 1,\ldots,q$, where $c$ is the generator of the center $Z(\tga) = \Ker(\tga\to \ga)$ with the remaining  defining relations
\begin{equation}\label{gen}
\tilde{g_1}^{e_1}= \ldots = \tilde{g}_n^{e_n } = c, \  \tilde{g}_1\cdots \tilde{g}_n [\tilde{a}_1,\tilde{b}_1]\cdots [\tilde{a}_q,\tilde{b}_q] = c^{2-2q-n}.
 \end{equation}
Under the canonical map $\tga\to \ga$, the generators $\tilde{g_i}, \tilde{a}_j,\tilde{b}_j$ are mapped to the generators $g_i, a_j,b_j$ (see \cite{LLR}, \cite{Scott}). We will be mostly interested in groups of signature $(0;e_1,\ldots,e_n)$. They admit a presentation as above with the elements $a_1,b_i$ absent.
This presentation makes sense also for finite subgroups of $\SU(2)$. In this case, 
$n\le 3$ and $e_1^{-1}+\ldots+e_n^{-1} >1$ and $c$ automatically satisfies $c^2 = 1$.

From now on we will assume that 
$$\bbD/\ga \cong \bbP^1$$
i.e. the signature  of $\ga$ is equal to $(0;e_1,\ldots,e_n ),$ where
$$\mu_\Gamma:= n-2-\sum_{i=1}^n\frac{1}{e_i} > 0. $$
Let $\ga' \subset \ga$ be a torsion-free normal subgroup of finite index (it always exists). The quotient $C = \bbD/\ga'$ is a compact Riemann surface of some genus $g > 1$ and $\Gamma' $ is isomorphic to the fundamental group of $C$. We will call $\Gamma'$ a \emph{surface subgroup of 
$\Gamma$ of genus $g$}.  The group  $G = \Gamma/\Gamma'$ acts faithfully on $C$ by holomorphic automorphisms.  Let 
$$V(C) = \Spec \bigoplus_{m=0}^\infty H^0(C, \calO_C(mK_C)).$$
If $C$ is not hyperelliptic, this is the affine cone over $C$ embedded canonically in $\bbP^{g-1}$. The group $G$ acts naturally on $V(C)$ and we obtain 
\beq
V_\Gamma \cong V(C)/G.
\eeq
This isomorphism allows one to find a minimal smooth compactification of the surface $V_\Gamma$. To avoid technicalities, we assume that $C$ is not hyperelliptic. First we replace $V(C)$ with the projective cone over $C\hookrightarrow \bbP^{g-1}$ and then take its natural minimal resolution. The resulting projective surface $\overline{V(C)}$ is isomorphic to the projectivization of the tangent line bundle $T_C$. It has the exceptional section $S_0$ with self-intersection $2-2g$ and the section at infinity $S_\infty$. The group $G$ acts on $\overline{V(C)}$ with finitely many fixed points: $n$ orbits of fixed points $x_1^0,\ldots,x_n^0$ lying on $S_0$ and $n$ orbits of fixed points $x_1^\infty,\ldots,x_n^\infty$ lying on $S_\infty$. Each pair $(x_i^0,x_i^\infty)$ lies on the same fibre of the projection $\overline{V(C)}\to C$. Choose a generator $g_i$ of $G_{x_i}$. Then $g_i$ has eigenvalues  $\eta_{e_i},\eta_{e_i}$ (resp.   $\eta_{e_i},\eta_{e_i}^{-1}$) at $x_i^0$ (resp. $x_i^\infty$), where $\eta_k$ denote a primitive $k$th root of the unity.

The quotient surface $\overline{V(C)}/G$ has $n$ singular cyclic quotient singularities  on each curve  $\bar{S}_0 = S_0/G$ and $\bar{S}_\infty = S_\infty/G$. Let
$$\pi:\overline{V_\Gamma } \to \overline{V(C)}/G$$ be a minimal normal crossing resolution of singularities. Let $\bfE_0$  and $\bfE_\infty$ be the pre-images of the curves $S_\infty/G$ and $S_0/G$. Then 
$$ \overline{V_\Gamma }\setminus \bfE_\infty \to (\overline{V(C)}\setminus S_\infty)/G = V_\Gamma $$
is a minimal normal crossing resolution of the affine surface $V_\Gamma $ with the exceptional curve equal to $\bfE_0$. We have 
$$  \overline{V}_\Gamma = V_\Gamma^*\cup \bfE_\infty \cup \bfE_0.$$
The surface $\overline{V_\Gamma }$ does not depend on the choice of $\ga'$. It a unique minimal normal crossing smooth $\bbC^*$-equivariant compactification of $V(\ga)^*$.

Each of the curves $\bfE_0,\bfE_\infty$ is the union of smooth rational curves intersecting each other transversally at  most one point with star-shaped incidence graph with $n$ arms and a central vertex (we draw only the case  $n = 3$). The curve $\bfE_0$ looks as in the next picture.

\begin{figure}[ht]
\xy@={(0,0),(13,0),(26,0),(13,-13)}@@{*{\bullet}};
(0,0)*{};(26,0)**{}**\dir{-};
(13,0)*{};(13,-13)**{}**\dir{-};
(0,3)*{-e_1};(13,3)*{2-n};(26,3)*{-e_2};(17,-13)*{-e_3};
(-50,0)*{};(0,10)*{};
\endxy
\caption{$\bfE_0$}\label{d}
\end{figure}

The numbers indicate the self-intersection of curves. 

\begin{figure}[ht]
\xy@={(0,0),(10,0),(30,0),(40,0),(50,0),(70,0),(80,0),(40,-10),(40,-30)}@@{*{\bullet}};
(0,0)*{};(16,0)**{}**\dir{-};
(24,0)*{};(56,0)**{}**\dir{-};
(66,0)*{};(80,0)**{}**\dir{-};
(40,0)*{};(40,-16)**{}**\dir{-};
(40,-24)*{};(40,-30)**{}**\dir{-};
(20,0)*{\cdots};(60,0)*{\cdots}; (40,-18)*{\vdots};
(-20,0)*{};(0,10)*{};
\endxy
\caption{$\bfE_\infty$}\label{dd}
\end{figure}

Each  arm consists of $e_i$ smooth rational curves, all curves have   self-intersection equal to $-2$.

Note that in the case when $\Gamma$ is a polyhedral group, the exceptional curve of the minimal resolution looks like $\bfE_\infty$ and the graph is a Dynkin diagram. The curve at the infinity looks like $\bfE_0$.

 The projection $\overline{V(C)}\to C$ induces a structure of a non-minimal ruled surface on $\overline{V}_\Gamma$:
$$f: \overline{V}_\Gamma = \overline{V(C)}/G \to C/G = \bbP^1.$$
Let $p_1,\ldots,p_n\in \bbP^1$ be the exceptional orbits of $\Gamma$ on $\bbD$ with stabilizer subgroups  of orders $e_1,\ldots,e_n $. The fibre over $p_i$ in $\overline{V}_\ga$ is a chain of curves with  incidence graph

\medskip
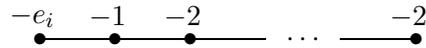
\xy @={(0,0),(10,0),(20,0),(50,0)}@@{*{\bullet}};
(0,0)*{};(30,0)**{}**\dir{-};
(40,0)*{};(50,0)**{}**\dir{-};
(35,0)*{\cdots};
(-1,3)*{-e_i};(9,3)*{-1};(19,3)*{-2};(49,3)*{-2};
(-30,0)*{};(0,-5)*{};
\endxy Here the number of $(-2)$-curves is equal to $e_i-1$ (unless $e_i = 2$ in which case it is equal to 2). All other fibres are isomorphic to $\bbP^1$. The central components of the curves $\bfE_0$ and $\bfE_\infty$  are the disjoint special sections $s_0$ and $s_\infty$ of the ruling. 

Let $F_i$ be the $(-1)$-curve in the fibre over $p_i$. Blow down the curves $F_i$ and then  blow down the components of  $\bfE_\infty$ contained in the fibre. We  arrive at a minimal ruled surface $\bfF_{n-2}$. Its exceptional section is the image of $s_0$. In the case $n= 3$, we can blow down the exceptional section to get a birational morphism 
$$\sigma:\overline{V}_\Gamma\to \bbP^2.$$
The image of the section $s_0$  is a point $p_0$. The images of the fibres $f^{-1}(p_i)$ are  members $\ell_1,\ell_2,\ell_3$ of the pencil of lines through $p_0$. The image of the section $s_\infty$ is a line $\ell$ not containing $p_0$. The surface $\overline{V}_\Gamma$ is obtained by a sequence of the blow-ups of the points $\ell_i\cap \ell$ and their infinitely near points.

\section{The dual McKay correspondence}
Recall that the ``easy part'' of the McKay correspondence for binary polyhedral groups is a natural bijection between the set $\calC(\Pi)$ of  conjugacy classes of non-trivial elements of the binary polyhedral group $\tilde{\Gamma}$ and the irreducible components of the curve $\bfE_0$. There are many ways to establish this correspondence. For example, we can use the description of the generators of $\Gamma$ in terms of natural generators of the fundamental group of  the boundary $\calT$ of a tubular neighborhood  of $\bfE_0$ on the surface $\overline{V}_\ga$ given by D. Mumford \cite{Mu}. The boundary  $\calT$ is mapped to $\bfE_0$ and its restriction over each component $E$ of $\bfE_0\setminus \text{singular points}$ is a circle bundle and the image of $\pi_1(\text{fibre})$ in $\pi_1(\calT) \cong \pi_1(V^*)$ is a certain element $s_E$ of $\pi_1(\calT) \cong \pi_1(V_\Gamma^*)$. Mumford proves that the element $s_E$ corresponding to the central component  can be taken for the generator $c$ of the center of $\Pi$. The element $s_E$ corresponding to  an end-component $E$ can be taken as generators $\gamma_i$. The element $e_E$ corresponding to the component next to an end-component $E_i$ is equal to $s_{E_i}^2$, and so on. Thus the bijective correspondence follows from the following.

\begin{lemma} [J.-L. Brylinski \cite{Br}]  Let $\Pi\subset \SL(2,\bbC)$ be a noncyclic binary polyhedral group generated by elements $\gamma_1,\ldots,\gamma_r,c$ with the standard defining relations  \eqref{gen}.
Any non-trivial non-central element is conjugate to a unique power $\gamma_i^t$ for some $i$ and $t< e_i$.
\end{lemma}

\begin{proof} Considering the action of $\Pi$ on $\bbP^1$, it is easy to see that any element of finite order is conjugate to a power  $\gamma_i^t$ as in the assertion of the lemma. To prove the uniqueness, one considers the standard action of $\tilde{\Gamma}$ on $\bbC^2$ and compares the characters of each power of $\gamma_i$.
\end{proof}

To extend the previous lemma to the hyperbolic case we first modify the definitions by considering everything modulo the center. 

\begin{lemma}   Let $\Pi\subset \widetilde{\SU}(1,1)$ be a  subgroup generated by elements $\gamma_1,\ldots,\\
\gamma_r,c$ with standard  relations \eqref{gen}. Any non-trivial  element of finite order modulo the center is conjugate to a unique element  $\gamma_i^vc^u$ for some $1\le i\le r$ and $0<v< e_i$ and $u \ge 0$.
\end{lemma}

\begin{proof} Again considering the action of $\Pi$ on the upper-half plane it is easy to see that any element in $\Gamma = \Pi/\la c\ra$ is conjugate to a power of the coset of a unique  $\gamma_i$.  It remains to see that the elements $\gamma_i^vc^u$ and $\gamma_i^{v'}c^{u'}$ are not conjugate if $(v,u) \ne (v',u')$. This is obvious. First we may assume that $u = 1$. Raising the elements to the  $e_i$th power we obtain that $c^{v}$ is conjugate to $c^{v'+e_iu'}$. Since $v < e_i$ this implies that $u' = 1, v= v'$.
\end{proof}

It follows from the last lemma that we can identify the irreducible components of the curve with conjugacy classes of elements $\gamma_i^v, 0 < v <e_i$ and $c$ such that each element of finite order modulo the center is conjugate to one of those modulo the center. 

\section{$G$-linearized vector bundles}
Here we recall some known facts about $G$-linearized vector bundles. Let $X$ be a nonsingular projective curve and $E\to X$  an algebraic vector bundle of rank $r$ on $X$. Let $G'$ be a finite group acting on $X$ not necessarily faithfully. A $G'$-linearization of $E$ is a lift of the action of $G'$ on $X$ to an action on the total space of $E$. Let $G$ be the image of $G$ in $\Aut(X)$. Obviously, $E$ is $G$-invariant, i.e., for any $g\in G$ there is an isomorphism $\phi_g:g^*(E)\to E$. A  $G$-linearization on $E$   is a choice of compatible sets of isomorphisms $\phi_g$. It may not exist, but it always exists if we replace $G$ with some 
central extension $G'$. A well-known argument introduces the group
$\calG_E$ whose elements are the pairs $(g,\phi_g), g\in G, \phi_g:g^*(E)\overset{\sim}{\to} E$. The group law is 
$$(g,\phi_g)\cdot (g',\phi_{g'}) = (gg', \phi_{g'}\circ g'{}^*(\phi_g)).$$
The projection $(g,\phi_g) \mapsto g$ defines a surjective homomorphism $\calG_E\to G$. Its kernel is the group of automorphisms of $E$ identical on the zero section. For example, when $E$ is simple, the kernel is equal to $\bbC^*$. A $G'$-linearization of $E$  defines a homomorphism of $G$-extensions 
$G'\to \calG_E$. If $G'$ is a faithful action on $E$, then the homomorphism $G'\to \calG_E$ is injective. Conversely, any finite subgroup of $\calG_E$ that covers $G$ defines a faithful linearization of $E$.

Assume that $E$ is a simple vector bundle. Then we have a central extension
$$1\to \bbC^* \to \calG_E \to G \to  1.$$
Its class $\xi$ is an element of the finite group $M(G):= H^2(G,\bbC^*)$, the Schur multiplier of $G$. Let $m$ be the order of $\xi$. The Kummer exact sequence
$$1\to \mu_k \to \bbC^* \to \bbC^* \to 1$$
gives an exact sequence
$$0\to \Char(G)/m\Char(G) \to H^2(G,\mu_m) \to  H^2(G,\bbC^*)[m]\to 0,$$
where $\Char(G) = \Hom(G,\bbC^*)$ and $A[m]$ denotes the $m$-torsion part of an abelian group $A$. This shows that  $\xi$ is the image of an element of $H^2(G,\mu_m)$ that defines a non-trivial central extension
\beq\label{newext}
1\to \mu_m \to m.G \to G\to 1
\eeq
such that $\calG_E = \bbC^*\times_{\mu_m} m.G$. In particular, $m.G$ defines a linearization of $E$ and $E$ cannot be linearized by any smaller group. Note that the extension \eqref{newext} is not uniquely defined.

Let $E$ be a $G$-linearized rank $r$ vector bundle over a nonsingular projective curve $X$. For any point  $x$  in $X$ its stabilizer subgroup $G_x$ acts naturally on the fibre $E_x$ of $E$. The corresponding linear representation 
$G_x\to \GL(E_x)$ is called the \emph{isotropy action}. We identify $X$ with the zero section of the total space $\bbV(E)$ of $E$. The action of $G_{x}$ on the tangent space $T_{x}(\bbV(E))$  decomposes into the direct sum of two invariant subspaces $T_x(X)$ and $T_x(E_x) = E_x$. Let $e(x) = |\bar{G}_{x}|$, where $\bar{G}_x$ is the image of $G_x$ in $\Aut(X)$. Let $\xi_{e(x)}$ be a primitive $e(x)$th root of unity. We may assume that the eigenvalue of a generator $g_x$ of $\bar{G}_x$ on $T_x(X)$ is equal to $\xi_{e(x)}$ and the eigenvalues of $g_x$ on $T_x(E_x)$ are equal to $\xi_{e(x)}^{q_i(x))}$ for some $0\le q_1(x)\le \ldots\le q_r(x) < e$.  

Let $Gx_1,\ldots,Gx_n$ be the orbits of $G$ lying over the branch points $p_1,\ldots,p_n$ of the projection $\pi:X\to Y = X/G$. Let $e_j = e(x_j), q_i^{(j)} = q_i(x_j), 1\le i \le r, 1\le j\le n$.  We call the vectors
$$\bfs_j = (q_1^{(j)}/e_j,\ldots,q_r^{(j)}/e_j), \ j = 1,\ldots, n.$$
the \emph{Seifert invariants} of $E$. Of course, an  expert would  rephrase this in terms of the degree of the corresponding parabolic bundle on $\bbP^1$ but we skip it.

Let $\pi_*(E)^G$ be the sheaf on $Y$ whose sections on an open subset $U$ are $G$-invariant sections of $E$ over $\pi^{-1}(U)$. Since $Y$ is a nonsingular curve, the sheaf $\pi_*(E)$ is locally free, hence $\pi_*(E)^G$ is locally free. Thus $\pi_*(E)^G$ is the sheaf of sections of a vector bundle  on $Y$. We denote it by $E^G$. Note that its total space is not in general isomorphic to  $\bbV(E)/G$. However it is isomorphic to this quotient over the subset of $Y$ over which the quotient map is unramified. In particular the rank of $E^G$ is equal to the rank of $E$.

Let $U$ be the universal cover of $X$ and $\pi_1(X,x)$ be the fundamental group of $X$ and   $\pi_1(X;G)$ be   given by the extension
\beq\label{fg}
1\to \pi_1(X,x) \to \pi_1(X;G) \to G\to 1.
\eeq
It is a subgroup of $G\times \Aut(U)$ which consists of pairs $(g,\tilde{g})$, where $g\in G$ and $\tilde{g}$ is a  lift of $g$ to an automorphism of $U$.
  
 \begin{remark} To follow the modern trend one can translate everything from above into the language of stacks or orbifolds. A $G$-linearized vector bundle on $X$ is a vector bundle on the smooth 1-dimensional DM-stack $[X/G]$. It is a reduced stack if $G$ acts faithfully.  The quotient variety $X/G$ is the coarse moduli space of the stack $[X/G]$. The group $\pi_1(X;G)$ defined in \eqref{fg} is the fundamental group of the stack $[X/G]$ (see \cite{BN}).  Any smooth projective DM-curve  is defined, up to isomorphism, by a smooth projective curve $Y$ of genus $q$, a (maybe empty) collection of points $y_1,\ldots,y_n$ on $Y$, a collection of integers $e_1,\ldots,e_n > 1$ and the generic inertia group $H$. There are three types of such DM-curves:hyperbolic, euclidean and spherical.  We will be interested only in hyperbolic DM-curves. A hyperbolic DM-curve is characterized by the condition 
 $$\chi(Y): = 2-2q-n +\sum_{i=1}^ne_i^{-1} <  0.$$
In the category of analytic stacks a hyperbolic DM-curve is isomorphic to the quotient stack $[C/H.G]$, where $C$ is a compact Riemann surface of genus $g > 1$ and $H.G$ is a finite extension of a  group $G$ of automorphisms of $C$  such that the projection $C\to C/G = Y$ is ramified over $p_1,\ldots,p_n$ with ramification indices $e_1,\ldots,e_n$.  There is a  quotient stack 
$[\bbD/H.\Gamma]$, where $H.\Gamma$ is an extension of a cocompact subgroup $\Gamma$ of $\PSU(1,1)$ of signature $(q;e_1,\ldots,e_n)$. It  is the universal cover of the DM-stack 
$[C/H.G]$, where $C = \bbD/\Gamma'$ for some surface subgroup of $\Gamma$ of genus $g$. In this paper we will be dealing with the case $q = 0$ and $H.G$ a central extension by a finite cyclic group $H$.
\end{remark}

Let  $\rho: \pi_1(X;G)\to \GL(r)$ be a linear representation. The restriction of $\rho$ to the subgroup $\pi_1(X)$ defines a vector bundle $E_\rho$ on $X$ equal to the quotient of the trivial bundle $U\times \bbC^r$ by the action $\gamma:(u,v) \mapsto (\gamma(u), \rho(\gamma)(v))$. The factor group $G = \pi_1(X;G)/\pi_1(X)$ acts naturally on $E_\rho$ and gives it the structure of a $G$-linearized vector bundle on $X$. The following result of A. Weil (see \cite{Gr}) gives a necessary and a sufficient condition for a $G$-linearized bundle to be isomorphic to a $G$-linearized bundle of the form $E_\rho$.

\begin{theorem}[A. Weil]  Let $E$ be a $G$-linearized vector bundle on $X$ and $E =\oplus _iE_i$ be its decomposition as a direct sum of indecomposable $G$-bundles.  Then $E \cong E_\rho$ for some linear representation  $\rho:\pi_1(X;G)\to \GL(r)$ if and only if for any $E_i$ the sum of the coordinates of its Seifert vectors is equal to the minus the degree of the vector bundle $E_i^G$.
\end{theorem}

It follows from this theorem that  $E_\rho$ is indecomposable if and only if $\rho$ is an irreducible representation.

\begin{definition} Let $G$ be a group acting on $X$ not necessarily faithfully. A $G$-linearized vector bundle $E$ is called $G$-stable (resp. $G$-semi-stable) if for any $G$-linearized subbundle $F$  of $E$ of smaller rank one has $\mu(F) < \mu(E)$ (resp. $\mu(F) \le \mu(E)$), where as usual the slope $\mu$ is defined as the ratio $\deg (E)/\rank (E)$.  If the equality takes place we say that $E$ is $G$-semistable.
\end{definition}

A $G$-stable $G$-linearized vector  bundle is always semi-stable as a vector bundle, but  not necessarily stable. Conversely, a semi-stable vector bundle but not stable vector bundle can be stable as a $G$-linearized bundle.

A theorem of Narasimhan and Seshadri \cite{NS} asserts that any irreducible unitary representation $\rho$ defines a stable vector bundle $E_\rho$ of degree 0. Conversely, any such bundle arises from an irreducible unitary representation of $\Gamma'$. A reducible unitary representation defines a semi-stable bundle such that each indecomposable summand is a stable bundle of degree 0. This correspondence makes an equivalence of the corresponding categories.

It is easy to extend the theorem of Weil to stable vector bundles (see \cite{Do1}).

\begin{theorem}  Let $\rho:\pi_1(X;G) \to U(r)$ be an irreducible  unitary linear representation of $\pi_1(X;G)$. Then the $G$-linearized vector bundle $E_\rho$ is $G$-stable.  It is stable if the restriction of $\rho$ to $\pi_1(X)$ is irreducible. Conversely, a $G$-linearized $G$-stable bundle satisfying the condition on its Seifert vectors from the Weil theorem is isomorphic to a bundle $E_\rho$ for some unitary irreducible representation of $\pi_1(X;G)$. It is stable if the restriction of the representation to $\pi_1(X)$ is irreducible.
\end{theorem}

In the next section we will interpret the condition from the Weil Theorem as the vanishing of the first Chern class of a $G$-linearized bundle.

\section{The second cohomology group}
Most of the results in this section are known in the theory of Seifert 3-manifolds (see \cite{LLR}). We recall them in a way convenient for an algebraist. 

Let $E$ be a line bundle on a nonsingular projective curve $X$ with the sheaf of sections $\calO_X(D)$ for some divisor $D$ on $X$. A $G$-linearization on $E$ is a choice of a divisor $D'$ in the linear equivalence class of $D$ such that $g^*(D') = D'$. Two such divisors define isomorphic $G$-linearized line bundles if and only if they differ by the divisor of a $G$-invariant rational function on $X$. The isomorphism classes of $G$-linearized line bundles form the group $\Pic(X;G)$ (the Picard group of the DM-stack $[X/G]$) which fits in the exact sequence
\beq
0\to \Char(G) \to \Pic(X;G) \to \Pic(X)^G \to H^2(G,\bbC^*) \to 1
\eeq
The image of the map $\Char(G) = \Hom(G,\bbC^*) \to \Pic(X;G)$ consists of trivial line bundles with the action defined by a character of $G$. The last non-trivial map defines an obstruction class for the existence of linearization on a $G$-invariant line bundle.

Also we have  an exact sequence
\beq\label{exa1}
0\to \Pic(Y) \to \Pic(X;G) \to \bigoplus_{j=1}^n\Char(G_{x_i}) \to 0,
\eeq
where  the first non-trivial map is defined by taking the pre-image of a line bundle under the projection $\pi:X\to Y$. The next  map is defined by the isotropy action.

 We identify the degree homomorphism $\deg:\Pic(Y)\to \Z$ with the first Chern class of a line bundle taking for the generator of $H^2(Y,\Z)$ the fundamental class of $Y$. One can define the first Chern class of a $G$-linearized bundle that gives a homomorphism $c_1^G:\Pic(X;G) \to H^2(\pi_1(X;G),\Z)$.
When $G$ acts faithfully on $X$, the group $\Gamma$ coincides with the orbifold fundamental group of $Y$. One has an exact sequence similar to \eqref{exa1}
\beq\label{exa2}
0\to H^2(Y,\Z) \to H^2(\pi_1(X;G),\Z) \to \bigoplus_{j=1}^nH^2(G_{x_i},\Z) \to 0,
\eeq
where the second map is defined by the restriction homomorphism for cohomology groups. The Chern class defines the commutative diagram 

\xymatrix{ &0\ar[d]& 0\ar[d]&& \\
0\ar[r]& \Pic^0(Y)\ar[r]\ar[d]& \Pic^0(X;G)\ar[r]\ar[d]&0\ar[d]& \\
0\ar[r]& \Pic(Y)\ar[r]\ar[d]_{c_1}& \Pic(X;G)\ar[r]\ar[d]_{c_1^G}&\bigoplus_{j=1}^n\Char(G_{x_i})\ar[r]\ar[d]_{\cong}& 0\\
0\ar[r]& H^2(Y,\Z)\ar[r]\ar[d]& H^2(\pi_1(X;G),\Z)\ar[d]\ar[r]&\bigoplus_{j=1}^nH^2(G_{x_i},\Z)\ar[r]\ar[d]& 0\\ 
& 0& 0&0& }
Here the isomorphism $\Char(G_{x_i})\to H^2(G_{x_i},\Z)$ is defined by using the Euler exact sequence 
$$0 \to \Z \to \bbC \overset{\exp}{\longrightarrow} \bbC^* \to 0.$$

In the case relevant to this paper, $Y = \bbP^1$ and hence $\Pic^0(Y) = 0$. Thus the Chern classes maps are isomorphisms. The group $H^2(\pi_1(X;G),\Z)$ is easy to compute. We have 
$$H^2(\pi_1(X;G),\Z) \cong \Z\oplus \Tors(H^2(\pi_1(X;G),\Z)).$$
To compute the torsion group we use again the Euler sequence. It gives 
\[\small{\Tors(H^2(\pi_1(X;G),\Z)) \cong \Tors(H^1(\pi_1(X;G),\bbC^*)) = \Tors(\Hom(\Gamma/[\Gamma,\Gamma],\bbC^*)).}\]
Assume that $X = \bbD/\Gamma'$, where $\Gamma'$ is a surface subgroup of $\Gamma$ and $G = \Gamma/\Gamma'$. In this case $G$ acts faithfully on $C$ and  $\pi_1(X;G) = \Gamma.$  Using the presentation of $\Gamma$ we find that 
$$\Gamma/[\Gamma,\Gamma] = \Coker(\Z^{n+1}\overset{f}{\to} \Z^{n }),$$
where the map $f$ is given by the matrix
$$\begin{pmatrix}1&e_1&0&\ldots&&0\\
1&0&e_2&0&\ldots&0\\
\vdots&\vdots&\vdots&\vdots&\vdots\\
1&0&0&0&\ldots&e_n \end{pmatrix}.
$$
Using the elementary theory of abelian groups we obtain
$$\Gamma/[\Gamma,\Gamma] \cong   \Z/a_1\Z\oplus \ldots \Z/a_n\Z,$$
where $a_i = c_i/c_{i-1}$,  $c_0 = 1$, and
$$c_k = g.c.d. ((e_{i_1}\cdots e_{i_k})_{1\le i_1 <\cdots < i_k\le n}), \ k = 1,\ldots, n-1.$$
In particular, this group is finite and coincides with its torsion group. We have
\beq\label{gammaab}
|\Tors(H^2(\Gamma,\Z))| = |\Gamma/[\Gamma,\Gamma] | = c_{n-1}.
\eeq
It is trivial if the numbers $e_1,\ldots,e_n $ are coprime. 

Let us identify $H^2(\bbD/\Gamma,\Z)$ with $\Z$ such that $1$ is the fundamental class of $\bbD/\Gamma$. Let $\gamma_0$ be a generator of the free part of $H^2(\Gamma,\Z)$ such that
$$\alpha(1) = l\gamma_0, \ l > 0.$$
Then
$$\Z/l\Z\oplus \Tors(\Gamma/[\Gamma,\Gamma])\cong \bigoplus_{k=1}^n\Z/e_k\Z.$$
In particular,
\beq\label{d}
l = \frac{e_1\cdots e_n }{c_{n-1}} = \lcm (e_1,\ldots,e_n ).
\eeq
This suggests to identify the free part of $H^2(\Gamma,\Z)$ with the subgroup $\frac{1}{l}\Z$ of $\bbQ$, so that the homomorphism $\alpha$ becomes the natural inclusion $\Z\hookrightarrow \frac{1}{l}\Z$
(cf. \cite{LLR}, Lemma 4.2).

Consider the isomorphism 
$$c_1^G:\Pic(X;G) \to H^2(\Gamma,\Z).$$
Recall that for any element $L$ of $\Pic(X;G)$ we can assign the degree of the vector bundle $L^G$ on $X/G$ and the Seifert numbers $q_j/e_j, j = 1,\ldots,n$. Let $D_j$ be the reduced divisor with support equal to the fibre over the branch point $p_j$ of $p:X\to Y=  \bbP^1$. Obviously it is a $G$-invariant divisor and hence $L_j = \calO_X(D_i)$ is a $G$-linearized line bundle. It follows from the Hurwitz formula that we have an isomorphism of $G$-linearized line bundles
\beq\label{old}
K_X \cong p^*(K_{\bbP^1})\otimes \bigotimes_{j=1}^nL_j^{e_j-1}.\eeq
The Seifert numbers of $K_X$ are equal to $(\frac{e_i-1}{e_1},\ldots,\frac{e_n-1}{e_n})$. Thus the image of  $c_1(L_j)$ in the  group $\oplus_{j=1}^ne_i^{-1}\Z/\Z$ is equal to $(0,\ldots,0,1/e_j,0,\ldots,0)$.  This shows that any $G$-linearized line bundle is isomorphic to the bundle
$$L \cong p^*(\calO_{\bbP^1}(m))\otimes \bigotimes_{j=1}^nL_j^{q_j},$$
where $q_j/e_j$ are the Seifert numbers of $L$ and $m$ is an integer. Since $L_j^{q_j}$ is a subsheaf of $\calO_X(e_jD_j) \cong p^*\calO_{\bbP^1}(1)$ and also contains $\calO_X$ we easily see that 
$(L_j^{q_j})^G$ is a subsheaf of $ \calO_{\bbP^1}(1)$ and contains $\calO_{\bbP^1}$. This easily implies that $(L_j^{q_j})^G \cong \calO_{\bbP^1}$ (see \cite{EL}, p.39). The projection formula shows that 
$$L^G \cong \calO_{\bbP^1}(m).$$
Comparing the degrees we find
$$\deg L = |G|m+\sum_{j=1}^nq_j|G|/e_j = |G|(m+\sum_{j=1}^n\frac{q_j}{e_j}).$$
The number
$$\vdeg(L):= m+ \sum_{j=1}^n\frac{q_j}{e_j} = \deg L/|G|$$
is called the \emph{Chern number} (or the \emph{virtual degree}) of a $G$-linearized line bundle $L$.
For example,
$$\vdeg (K_X) = \mu(\Gamma):=  n-2-\sum_{j=1}^n\frac{1}{e_j}.$$
($= \frac{1}{2\pi}\times$ the area of the polygon with angles $\pi/e_i$ on the hyperbolic plane).

For any $G$-linearized vector bundle $E$ its determinant $\det E$ is a $G$-linearized line bundle. We set
$$\vdeg(E) = \vdeg(\det E).$$
Grothendieck shows in \cite{Gr} that 
$$\vdeg(E) = \deg(E^G) +\sum_{j=1}^n\sum_{i=1}^r\frac{q_j^{(i)}}{e_j}.$$
Thus we obtain that the condition for an indecomposable $E$ to be isomorphic to $E_\rho$ is that $\deg (\det E) = |G|\vdeg(\det E)= 0$. Note that $\det (E)^G \ne \det (E^G)$ since the sum of the coordinates of the corresponding Seifert vector $(q_j^{(1)}/e_j,\ldots,q_j^{(r)}/e_j)$ differs from the Seifert invariant of $E$ at $p_j$ by an integer.

Let $\gamma_0$ be a generator of the free part of $\Pic(X;G)$. The forgetful homomorphism $\Pic(X;G) \to \Pic(X)^G$ composed with the degree map corresponds to the restriction map $H^2(\Gamma,\Z) \to H^2(\Gamma',\Z) = H^2(X,\Z)$. It sends $K_X$ to $2g-2$. On the other hand, the composition  $\Pic(Y) \to \Pic(X:G)\to \Pic(X)$ is the multiplication by $|G|$.  This implies that 
$$c_1^G(K_X) - s\gamma_0\in \Tors(H^2(\ga,\Z)).$$
where
\beq\label{s}
s = l\frac{2g-2}{|G|} = \lcm(e_1,\ldots,e_n)(n-2-\sum_{i=1}^n\frac{1}{e_i}).
\eeq
In fact, one can show that a generator $\gamma_0$ can be changed in such a way that we have the equality
$$c_1^G(K_X) = s\gamma_0.$$
The generator $\gamma_0$ is equal to $c_1^G(L_0)$, where $L_0$ is $G$-linearized line bundle with $\vdeg(L_0) = \mu_\ga/s$ and the Seifert invariants $\frac{q_i}{e_i}$, where $sq_i\equiv 1\mod e_i$ (see \cite{Do3} or \cite{LLR}, Corollary 4.4). In the case when the surfcae $V_\ga$ is a hypersurface in $\bbC^3$ with equation $f(x,y,z) = 0$, where $f(x,y,z)$ is a quasi-homogeneous polynomial of degree $d$ with positive integer weights $w_1,w_2,w_3$ we have (see \cite{Anna},\cite{Do2}, \cite{Sch})
\beq\label{my}
s = d-w_1-w_2-w_3.
\eeq

Let $\Gamma'$ be a surface subgroup of $\ga$ of genus $g$ with $G = \ga/\ga'$. We apply the previous discussion by taking $C = X$. The class $\gamma$ of the extension \eqref{univext} defining our group $\Pi$ is equal to $c_1^G(K_C)$.  Recall that the natural homomorphism of the cohomology groups 
$H^2(\Gamma,\Z)\to H^2(\Gamma,\Z/m\Z)$ is defined in terms of the extension classes as follows. Given an extension
$1\to \Z \to A \to \Gamma\to 1$ with the class $e$, its image in $H^2(\Gamma,\Z/m\Z)$ is the class of the extension $1\to \Z/m\Z\to A/mA \to \Gamma\to 1$. The exact sequence 
$$0\to \Z \to \Z \to  \Z/m\Z\to 0$$ 
together with the fact that $H^3(\Gamma,\Z) = 0$ (\cite{LLR}, p. 3788) gives an isomorphism 
\beq\label{h3}
 H^2(\Gamma,\Z)/mH^2(\Gamma,\Z)\cong H^2(\Gamma,\Z/m\Z).\eeq
This  shows that the image of the class $e$ is equal to zero in  $H^2(\Gamma,\Z/m\Z)$ if and only if 
$e$ is divisible by $m$ in $H^2(\Gamma,\Z)$. In particular, we obtain that $\gamma = s\gamma_0$ is mapped to zero in $H^2(\Gamma,\Z/m\Z)$ if and only if $m$ divides $s$. 

On the hand, let us consider the image of $\gamma$ under the composition of 
$H^2(\Gamma,\Z)\to H^2(\Gamma,\Z/m\Z)$ and the restriction map $H^2(\Gamma,\Z/m\Z)\to H^2(\Gamma'/,\Z/m\Z) = \Z/m\Z$. The image of $\gamma$ is equal to $2g-2 \mod m$. So the extension 
\beq\label{splitt}
1\to \la c\ra/\la c^m\ra \to \Pi'/\la c^m\ra \to \Gamma' \to 1
\eeq
splits if and only if $m$ divides $2g-2$.

Recall that the group $\Pi'$ is generated by $c, a_i,b_i, i = 1,\ldots,g$ with defining relation
$[a_1,b_1]\cdots [a_g,b_g] = c^{2-2g}$. If $m|2-2g$ we can define an explicit distinguished splitting by sending the standard generators $x_i,y_i$ of $\ga' =\pi_1(C)$ to $a_i \mod \la c^m\ra, b_i\mod \la c^m\ra$.  We call such a splitting \emph{canonical splitting}. Its image is a subgroup of $\Pi'/\la c^m\ra$ which is normal in $\Pi/\la c^m\ra$.  The quotient group is a cyclic central extension $m.G$ of $G$. The extension 
\beq\label{splitt2}
1\to \la c\ra/\la c^m\ra \to \Pi/\la c^m\ra \to \Gamma \to 1
\eeq
can be reconstructed from $m.G$ by using an isomorphism 
\beq
 \Pi/\la c^m\ra \cong m.G\times_G \ga,
\eeq
where $\ga\to G$ is the quotient map $\ga \to \ga/\ga'$ and $m.G\to G$ is the defined by the extension. In other words, the class of the extension \eqref{splitt2} is the image of the extension $m.G$ under the inflation map $H^2(G,\Z/m\Z) \to H^2(\ga,\Z/m\Z)$.

\section{Admissible unitary representations of $\Pi$}\label{sec5}
As we have explained in the Introduction, the first side of the McKay correspondence must be a class of linear representations of  $\Pi = \pi_1(V_\Gamma^*)$. Since any representation of a finite group is conjugate to a unitary representation, it is natural to deal with unitary representations of $\Pi$. Any unitary representation of $\Gamma$ defines via the composition with $\Pi\to \Gamma$ a unitary representation of $\Pi$. 

The proof of the next lemma was kindly explained to me by Daniel Allcock.  

\begin{lemma} Let $\rho: \Pi\to \GL(r)$ be an irreducible linear representation. Then $\rho(c)$ is of finite order. There exists a surface subgroup $\Gamma'$ of $\ga$ of genus $g$ such that 
$$\rho(c)^{2g-2} = 1.$$
\end{lemma}

\begin{proof} It follows from the presentation \eqref{gen} of  $\Pi$ 
that $c^a \in [\Pi,\Pi]$ for some $a > 0$. Since the commutator of $\GL(r)$ is equal to $\SL(r)$ we obtain that $\rho(c)^a\in \SL(r)$. Since $\rho$ is irreducible, $\rho(c)^a$ belongs to the center of $\SL(r)$, hence  $\rho(c)^{ar} = 1$. This proves that $\rho(c)$ is of finite order.

Let $\ga_1$ be any surface subgroup of $\ga$ of some genus $g'$. Since $\ga_1/[\ga_1,\ga_1] \cong \bbZ^{2g'}$, there exists a subgroup $\ga_2$  of $\gamma_1$ of index divisible by the order $n$ of $\rho(c)$. The normal subgroup $\ga' = \cap_{g\in G} g\ga_2 g^{-1}$ of $\ga$ is a surface subgroup of $\ga$ of index divisible by $n$. By Hurwitz's formula, $n$ divides $2g-2$, where $g$ is the genus of $\ga'$. 

\end{proof} 

Since any unitary representation decomposes into the direct sum of irreducible representations, it follows from the previous lemma that the image of the center of $\Pi$ is always a finite cyclic group.

\begin{definition} A  unitary  linear representation $\rho:\Pi \to \U(r)$ is called {\it admissible} if its irreducible summands remain irreducible after restriction to any subgroup of finite index.  The order of the image of the center of $\Pi$ will be called  the \emph{level} of the representation.  Two admissible representation $\rho$ and $\rho'$ are called \emph{equivalent} if there exists a subgroup $\Pi'$ of $\Pi$ of finite index such that the restriction of $\rho$ and $\rho'$ to $\Pi'$ are  isomorphic.
\end{definition}

It is known that an irreducible representation is admissible if and only if  its image in $\GL(r)$ is Zariski-dense.

 \begin{remark} A linear representation of $\Pi$ that factors through an irreducible representation of a finite group is an irreducible unitary representation of $\Pi$ with $\rho(c)$ of finite order.  The reason why we would like to ignore such representation is that there are ``too many'' of these representations. For example, it is known that the group $\Gamma$ of signature $(0;2a,2b,c)$ surjects to all but finitely many symmetric groups $S_n$ \cite{Liebeck}.

\end{remark}

\begin{lemma}\label{prev}Let $\rho$ be an admissible  representation of $\Pi$ of level $m$. Assume that  the extension
$$1\to \la c\ra /\la c^k\ra \to \Pi/\la c^k\ra \to \Gamma \to 1$$
splits for some $k|m$. Then $\rho$ is equivalent to an admissible representation of level dividing $m/k$.
\end{lemma}

\begin{proof} Let $r:\Gamma \to  \Pi/\la c^n\ra$ be a splitting and $\widetilde{\Gamma}$ be the pre-image of $r(\Gamma)$ in $\Pi/\la c^m\ra$. It defines an extension
$$0 \to \Z/(m/k)\Z \to \widetilde{\Gamma} \to \Gamma \to 1.$$
This extension is isomorphic to the extension
$$0 \to \la c \ra /\la c^{m/k}\ra \to \Pi/\la c^{m/k}\ra \to \Gamma \to 1.$$
The  representation $\rho$ factors through $\Pi/\la c^m\ra$ and then restricts to $\widetilde{\Gamma}$ to define a representation of $\widetilde{\Gamma}$. Its composition with the projection $\Pi\to  \Pi/\la c^{m/k}\ra$ defines an admissible representation $\rho':\Pi\to U(r)$ of level dividing $m/k$. Let us show that $\rho$ and $\rho'$ are equivalent.

We have $\Pi/\la c^m\ra \cong \widetilde{\Gamma} \times \la c^k\ra /\la c^m\ra$. Let $\widetilde{\Gamma}'$ be the pre-image of $\widetilde{\Gamma}$ in $\Pi$. It is a subgroup of  finite index of  $\Pi$ that fits in  the extension
$$1\to \la c^m\ra \to \widetilde{\Gamma}' \to \Gamma \to 1.$$
It is clear that $\rho$ and $\rho'$ define the same representation of $ \widetilde{\Gamma}'$.
\end{proof}

\begin{proposition} Assume that $s$ ($=\lcm(e_1,\ldots,e_n)\mu(\Gamma)$)  is coprime to the order of $\Tors(H^2(\Gamma,\Z))$ ($= \frac{e_1\cdots e_n}{\lcm(e_1,\ldots,e_n)})$. Every  admissible representation of dimension $r$ is equivalent to an admissible representation of level $m$ dividing $r$ and coprime to $s$. 
\end{proposition}

\begin{proof} Since $\rho$ is irreducible and $c$ belongs to the center of $\Pi$, the image of $c$ is a diagonal matrix $[\alpha,\ldots,\alpha]$, where $\alpha^m = 1$. Let $f:\Pi\to U(1)$ be the composition of $\rho$ and the determinant map $U(r)\to U(1)$. 
It factors through $\Pi^{\ab} = \Pi/[\Pi,\Pi]$ and $f(c) = \alpha^{r}.$ The group $\Pi^{ab}$ is easy to compute since we know its presentation. It is equal to the abelian group defined by the square matrix 
$$\begin{pmatrix}1&e_1&0&\ldots&&0\\
1&0&e_2&0&\ldots&0\\
\vdots&\vdots&\vdots&\vdots&\vdots\\
1&0&0&0&\ldots&e_n \\
1&1&1&1&\ldots&1\end{pmatrix}.
$$
Its order is equal to the determinant of the matrix and we find
$$|\Pi^{\ab}| = e_1\cdots e_n (n-2-\sum_{i=1}^n \frac{1}{e_i}).$$
Note that this number coincides with the determinant of the intersection matrix of the exceptional curve $\bfE_0$. It follows from \cite{Mu} that the group $\Pi^{ab}$ is isomorphic to the group of connected components of the local Picard group of the singularity $(V_\ga,\{0\})$ (see \cite{Mu}).

Comparing this with the order of $\Gamma^{\ab}$ given in \eqref{gammaab}, we obtain
$$|\Pi^{\ab}| = s|\Gamma^{\ab}|.$$
This shows that the kernel of the natural map $\Pi^{\ab}\to \Gamma^{\ab}$ induced by the surjection $\Pi\to \Gamma$ is isomorphic to $\Z/s\Z$. Thus the image $\bar{c}$ of $c$ in $\Pi^{\ab}$ satisfies $\bar{c}^s = 1$, hence $f(c)^s = (\alpha)^{rs} =1$.  Since $m$ is the order of $\alpha$ in $\bbC^*$, we get 
$m\vert rs$. It remains to show that $\rho$ is equivalent to an admissible representation of level $m$ coprime to $s$.

Assume  $k = (m,s) > 1$. Recall from the previous section that $\Pi$ fits in the extension defined by the class $\gamma = s\gamma_0$, where $\gamma_0$ projects to a generator  of  $H^2(\Gamma,\Z)/\Tors$. By the assumption $(s,|\Tors(H^2(\Gamma,\Z)|) = 1$, hence the multiplication by $k$ is an automorphism of $\Tors(H^2(\Gamma,\Z))$. Thus $\gamma$ is divisible by $k$ in $H^2(\ga,\Z)$. Applying \eqref{h3} (with $m$ replaced by $k$) we obtain  that  the extension
$$1\to \la c\ra /\la c^n\ra \to \Pi/\la c^n\ra \to \Gamma \to 1$$
splits. It remains to apply  Lemma \ref{prev}.
\end{proof} 

Note that the proof also shows that we may replace the assumption   on $s$ with the assumption that $m$ is coprime to $|\Tors(H^2(\Gamma,\Z))|$.

\begin{theorem}\label{const}  Let $\Gamma'$ be a surface subgroup of $\ga$ of genus $g$  with finite quotient $G$ and $C = \bbD/\ga'$. Assume $m$ divides $2g-2$. There is a natural bijective correspondence between  admissible  $r$-dimensional representations of $\Pi$ of level dividing $m$ and  $m.G$-linearized stable rank  $r$ vector bundles on $C$, where the image of the class of the extension $m.G$ under the inflation map $H^2(G,\Z/m\Z) \to H^2(\Gamma,\Z/m\Z)$ is equal to the class of  extension \eqref{splitt2}.\end{theorem}

\begin{proof} Let $\rho:\Pi\to U(r)$ be an admissible representation of level $m$. Fix a surface subgroup $\Gamma'$ of genus $g$ such and $m$ divides $2g-2$. As we explained at the end of section 5, the extension  $m.\ga$ defined by \eqref{splitt2} is equal to the image of an extension $m.G$ under the inflation map $H^2(G,\ Z/m\Z)\to H^2(\ga,\Z/m\Z)$. 
Now the representation $\rho:\Pi\to U(r)$ factors through $\Pi/\la c^m\ra = m.\Gamma$. Since $m.\Gamma$ is defined by the extension $m.G$, the group $m.\Gamma$ is isomorphic to the fundamental group of the DM-curve $[\bbD/m.G]$. This defines a $m.G$-linearized $m.G$-stable vector bundle over $C$. Conversely,  such a vector bundle $E$ defines a unitary irreducible representation of $m.\Gamma \cong \Pi/\la c^m\ra$ and hence an irreducible representation $\rho$ of $\Pi$ of level dividing $m$. Suppose the restriction of $\rho$ to some subgroup $H$ of finite index is reducible. Then we find a surface subgroup $\Gamma''$ of finite index in both $H$ and $\Gamma'$. The pull-back $E'$ of $E$ to the cover $\pi:C' = \bbD/\Gamma' \to C$ is defined by a reducible unitary representation, hence $E$ is not stable but semi-stable. Let $K = \Aut(C'/C)$ and $F$ be a unique maximal semi-stable subbundle of $E$ (see \cite{Se}, p. 15). It is a $K$-invariant subbundle of $E'$ of degree 0. The subbundle $F^K$ of $E^K \cong \calE$ is of degree 0 and rank equal to rank of $F$. Hence it defines a destabilizing subbundle of $E$, contradicting its stability. This proves that $\rho$ is an admissible representation of $\Pi$. It is easy to see that the construction establishes a bijective correspondence between the two sets.
\end{proof}

 \begin{remark} We know that  $|G|$ is divisible by $l = \lcm (e_1,\ldots,e_n)$. In fact, $|G|$ must be divisible by $2l$ if $l$ is even and the number of $e_i$'s such that $l/e_i$ is odd. This condition is sufficient for the existence of a not necessary normal torsion-free subgroup of finite index \cite{EEK}. It is not known what is the minimal index of a surface subgroup of $\Gamma$. However, if $n = 3$, accepted the extended Riemann hypothesis, in most cases the group $\Gamma'$ can be found such that $G = \Gamma/\Gamma'$ is a simple group isomorphic to $\PSL_2(\bbF_q)$ \cite{Ku}. Here, ``in most cases'' means that, fixing two of the numbers $e_i$'s, the values of the remaining $e_j$ for which this assertion is not true belongs to a set of density 0. 
 \end{remark}

\begin{example} Let $\Gamma$ be the triangle group with $(e_1,e_2,e_3) = (2,3,7)$. In this case $C$ can be chosen to be the Klein curve of genus 3 isomorphic to  the modular curve $X(7)$. The group $G = L_2(7)$ is the simple group of order 168, and $\SL(2,\bbF_7)$ is its non-split double extension. The order of $M(G)$ is equal to 2. The number $s$ is equal to 1. Thus the possible levels of an admissible representation of $\Pi$ are $1,2,4$. In particular, the dimension $r$ of an admissible representation of level $m > 1$ must be even. It is known (see \cite{Do1}) that there are exactly two non-isomorphic $\SL(2,\bbF_7)$-linearized stable bundles  on $C$ of rank 2. Each of them defines an admissible representation of level $2$. Also it is known there are exactly two non-isomorphic $\SL(2,\bbF_7)$-linearized stable bundles  on $C$ of rank 3 which arise from a $G$-linearized stable bundle.  Each of them defines an admissible representation of level $1$.  
\end{example}

\section{The McKay correspondence}
Let $\rho:\Pi\to U(r)$ be an admissible representation of $\Pi$ of level $m$ and $\ga'$ be a surface subgroup of genus $g$ with $m$ dividing $2-2g$. Let $E_\rho$ be a $m.G$-linearized stable vector bundle on $C = \bbD/\ga'$ constructed in Theorem \ref{const}. Let $f:\overline{V(C)}\to C$ be the projectivization of the tangent bundle of $C$. Let $V(C)^*$ be the complement of the sections $S_0$ and $S_\infty$. It is a $\bbC^*$-bundle over $C$. We know that the extension \eqref{splitt2} splits. The canonical splitting  defines a surjection 
$$\pi_1(\overline{V(C)}^*) = \Pi'\to \Pi'/\la c^m\ra = \ga'\times \Z/mZ \to \Z/m\Z.$$
Let $X^*\to V(C)^*$ be the corresponding cyclic cover. It extends to a cyclic cover $p:X\to \overline{V(C)}$ of degree $m$ ramified over $S_0\cup S_\infty$. The group $m.G$ acts on $X$ with the center $\Z/m\Z$ acting as the group of deck transformations. The surface $X$ is a minimal ruled surface over $C$ isomorphic to the projectivization of a $G$-invariant $m$th root of the tangent bundle of $C$ (it does not admit a $G$-linearization unless $m$ divides $s$). Let $\tilde{E}_\rho$ be the pull-back of $E_\rho$ on $X$. It is a $m.G$-linearized vector bundle. The restriction of the vector bundle $E_\rho' = (\tilde{E}_\rho)^{\Z/m\Z}$ to  $V(C)^*$ coincides with the pull-back of $E_\rho$ on $\overline{V(C)}^*$ but differs on the whole $\overline{V(C)}$. 

Now we have constructed an extension $E_\rho'$ of the pull-back of $E_\rho$ on $V(C)^*$ to the whole $\overline{V(C)}$. Let $\pi:\overline{V(C)}\setminus S_0\to (\overline{V(C)}\setminus S_0)/G$ be the projection to the quotient. We know that the quotient has only cyclic quotient singularities of types $A_{e_i-1}$ lying on the image of the curve $S_\infty$. The minimal resolution of the quotient is the surface $\widetilde{V}_\ga = \overline{V}_\Gamma\setminus \bfE_0$. Now we are in the situation of the ordinary McKay correspondence. The standard procedure is to take the coherent sheaf $\pi_*^G\calE_\rho'$, where $\calE_\rho'$ is the sheaf of sections of $E_\rho'$. Its pull-back on the resolution modulo its torsion subsheaf is a locally free sheaf $\tilde{E}_\rho$ on $\widetilde{V}_\ga $ (see \cite{GV}). Equivalently, one can resolve the singularities of the quotient via the $G$-Hilbert scheme $H = \Hilb_G(\overline{V(C)}\setminus S_0)$ (see \cite{Reid}). It comes with a universal $G$-equivariant family  $\alpha:\calZ\to H$ whose fibres are $0$-dimensional subschemes  $Z$ of $\overline{V(C)}\setminus S_0$ with $H^0(Z,\calO_Z)$ isomorphic to the regular representation of $G$. Let $\beta:\calZ\to \overline{V(C)}\setminus S_0$ be the natural projection. Then $\alpha_*^G(\beta^*\calE_\rho')$ is a locally free sheaf isomorphic to $\tilde{E}_\rho$.

Thus we have defined a correspondence 
\beq\label{final}
\{\text{admissible representations of $\Pi$}\} \overset{\rho\mapsto \tilde{E}_\rho}{\longrightarrow} \{\text{vector bundles on $\widetilde{V}_\ga\}$}.
\eeq

Let $K_0(\Pi)$ be the Grothendieck group based on the set of admissible representations of $\Pi$ and 
$K_0(\widetilde{V}_\ga)$ be the Grothendieck group based on vector bundles on the surface $\widetilde{V}_\ga$.

\begin{conjecture} The correspondence \ref{final} defines an isomorphism of groups
$$K_0(\Pi) \cong K_0(\widetilde{V}_\ga).$$
\end{conjecture} 

Note that the group $K_0(\widetilde{V}_\ga)$ is easy to compute. We know that the surface $\overline{V}_\ga$ is isomorphic to the blow-up of $\bbP^2$ at collinear points $p_1,\ldots,p_n$. The pre-image of the line containing these points is the union of the curve $\bfE_\infty$ and $n$ exceptional curves $F_i$ of the first kind contained in the fibres of the projection to $\bbP^1$. Since we know how $K_0$-groups change under  blow-ups we immediately obtain the following.

\begin{proposition} 
\begin{itemize} 
\item[(i)] The map 
$$\phi:K_0(\widetilde{V}_\ga) \to \Z\oplus \Pic(\widetilde{V}_\ga), \ E\mapsto (\rank E, c_1(E))$$
is an isomorphism of abelian groups.
\item[(ii)] The group $\Pic(\widetilde{V}_\ga)$ is a free abelian group generated by the classes $\omega_i$ which are dual to the divisor classes of irreducible components of $\bfE_\infty$ with respect to the intersection pairing $\Pic(\widetilde{V}_\ga)\times \Pic(\widetilde{V}_\ga)^\infty$, where $\Pic(\widetilde{V}_\ga)^\infty$ is the subgroup of divisors supported on $\bfE_\infty$.
\end{itemize}
\end{proposition}

Note that in the ordinary McKay correspondence we have an analogous statement for the minimal resolution of $\bbC^2/\Pi$ (Proposition 2 in \cite{GV}). Moreover,   the images in $K_0(\widetilde{V}_\Gamma)$ of the elements of $K_0(\Pi)$ corresponding to irreducible representations $\rho_i$ are equal to $(d_i,\omega_i)$, where $d_i = \dim \rho_i$. If we identify $\omega_i$ with the negatives of fundamental weights of the root system defined by the intersection matrix of $\bfE_0$, then the numbers $d_i$ are the coefficients of the maximal root.  In the hyperbolic  case I do not know what are the elements in $K_0(\Pi)$ whose images in 
$K_0(\widetilde{V}_\ga)$ are projected to the elements $\omega_i$ of $\Pic(\widetilde{V}_\ga)$.

\section{Weighted projective lines} Here we briefly comment on a certain relation of our work with the canonical algebras of C. Ringel and  weighted projective lines  of W. Geigle and H. Lenzing. There is enormous literature on this subject, so we refer only to  \cite{CB} and  \cite{LP}, where further references can be found. 

Let $\Gamma$ be a fuchsian group of signature $(0;e_1,\ldots,e_n)$. The algebra of automorphic forms 
$A([\ga,\ga])$ is generated by $n$ elements $f_1,\ldots,f_n$ with relations  of the form
$$f_1^{e_1}+a_if_i^{e_i}+f_{n}^{e_n} = 0, \ i = 2,\ldots,n-1,$$
where $a_2 = 1,a_i\in \bbC^*, i = 3,\ldots,n-1$. This result belongs to H. Poincar\'e who called the functions $f_1,\ldots,f_n$ the \emph{Halphen functions} (\cite{Poi}, p. 237, \cite{Do4}). This result implies easily that the  \'etale cover $U\to V_\ga^*$ corresponding to the commutator group $[\Pi,\Pi]$ is isomorphic to the punctured affine surface $U_\ga^* = U_\ga\setminus \{0\}$, where $U$ is an affine complete intersection surface in 
$\bbC^{n}$ given by the equations
$$z_1^{e_1}+a_iz_i^{e_i}+z_{n}^{e_n} = 0, \ i = 2,\ldots,n-1.$$
This result was proven by J. Milnor \cite{Mi} (n= 3) and W. Neumann \cite{Ne} ($n\ge 3$) independently of Poincar\'e's result. 

The ring $R$ of regular functions on the surface $U_\ga$ is a finitely algebra with action of the group $$(\lambda_1,\ldots,\lambda_n)\in \bbC^*{}^n: \lambda_1^{e_1}=\ldots = \lambda_n^{e_n}.$$
This is the diagonalizable algebraic group $D(M)$, where $M$ is an abelian group with presentation 
$\{g_1,\ldots,g_n: e_1g_1 = \ldots=e_ng_n\}$ and $D$ denotes the dual algebraic group. Using the theory of elementary divisors for abelian groups we find that 
$$M \cong H^2(\ga,\bbZ).$$
 The action of $D(M)$ is equivalent to the grading 
$$R = \bigoplus_{m\in M}R_m$$
of the algebra $R$ by the group $M$. A \emph{weighted projective line}  is the ring $R$ together with the above grading. In notation of section 5, let $\gamma = c_1^G(K_X)\in H^2(\ga,\bbZ)$. It is equal to $s\gamma_0$, where $\gamma_0$ generates the free part of $H^2(\ga,\bbZ)$. The quotient group $M/\gamma M$ is isomorphic to $\Pi^{\text{ab}} = \Pi/[\Pi,\Pi]$. The corresponding finite group $D(M/\gamma M)$ acts on the affine surface $U_\ga$ with  quotient  isomorphic to our surface $V_\ga$. 
The action is free on $U_\ga^*$. So, in this sense the surface $U_\ga^*$ is a maximal unramified abelian cover of the surface $V_\ga^*$.

Let $\bbX$ be a weighted projective line as above. Lenzing and others consider  the category 
$\coh~ \bbX$ of coherent sheaves on $\bbX$. By definition, it is the category of $M$-graded $R$-modules modulo the Serre subcategory of finite length modules. Let $R'$ be the same ring but with the grading defined by the free subgroup $\gamma \bbZ$ of $M$. For any $M$-graded $R$-module $P$, the coherent sheaf $\calF = P^\sim$ restricted to $U_\ga^*$ is a $\Pi^{\ab}$-sheaf that descents to a coherent sheaf on $V_\ga^* = U_\ga/\Pi^{\ab}$. It corresponds to $\bbZ$-graded $A(\Gamma)$-module 
$P^{\Pi^{\ab}}$. Admissible unitary representations of $\Pi$ correspond to coherent sheaves on $\bbX$ which define  locally free sheaves on $V_\ga^*$ isomorphic to the pre-image of a vector bundle on the DM-stack $\bbD/\Pi$ defined by an admissible unitary representation of $\Pi$.

The Grothendieck group $K_0(\coh~\bbX)$ was computed by Geigle and Lenzing. It is isomorphic to $K_0(\widetilde{V}_\ga)$.  The intersection matrix of the curve $\bfE_\infty$ is the negative of  a symmetric Cartan matrix. It defines a Kac-Moody Lie algebra $\frak{g}$ and associated loop algebra $\calL\frak{g}$ with  root lattice $K_0(\coh~\bbX)$. Crowley-Boevey proves that indecomposable coherent sheaves on $\bbX$ correspond to positive roots, a unique indecomposable sheaf for a real root and infinitely many for an imaginary root. It would be very interesting to determine which roots correspond to admissible unitary representations of $\Pi$. 

Note that the category $\coh~\bbX$ is derived equivalent to the category of modules over the canonical algebra $C(a,e_1,\ldots,e_n)$ of Ringel. The canonical algebra is the path algebra with relations of the following quiver obtained from  our graph describing the curve $\bfE_\infty$ by adding  one vertex joined to the extreme vertices of all arms of the following graph (we borrowed the picture from \cite{LP}). 
\newpage

\begin{figure}[ht]
\begin{center}

\label{ }
\xymatrix{&&&\bullet\ar[r]^{\alpha_{12}}&\bullet&\ldots&\bullet\ar[r]^{\alpha_{1e_1-1}}&\bullet\ar[ddr]^{\alpha_{1e_1}}&\\
&&&\bullet\ar[r]^{\alpha_{22}}&\bullet&\ldots&\bullet\ar[r]^{\alpha_{2e_2-1}}&\bullet\ar[dr]_{\alpha_{2e_2}}
&\\
&&0\ar[ur]_{\alpha_{21}}\ar[uur]^{\alpha_{11}}\ar[dr]^{\alpha_{n1}}&\vdots&\vdots&\vdots&\vdots&\vdots&\omega\\
&&&\bullet\ar[r]^{\alpha_{n1}}&\bullet&\ldots&\bullet\ar[r]^{\alpha_{ne_n-1}}&\bullet\ar[ur]^{\alpha_{ne_n}}&}
\end{center}
\caption{ }
\end{figure}

Here the path algebra is taken with the relations
$$a_i\alpha_{ie_i}\cdots \alpha_{i2}\alpha_{i1} +\alpha_{1e_1}\cdots \alpha_{12}\alpha_{11}+\alpha_{ne_n}\cdots \alpha_{n2}\alpha_{n1}  = 0, \ i = 2,\ldots, n-1.$$

\bibliographystyle{amsplain}

\end{document}